\documentclass[11pt,reqno,twoside]{article}

\usepackage{mysty_arxiv}

\usepackage[latin1]{inputenc}
\usepackage{tcolorbox}

\usepackage{fullpage}
\usepackage[top=2.35cm, bottom=2.35cm, left=2.35cm, right=2.35cm]{geometry}

\usepackage{tikz}
\usetikzlibrary{shapes.geometric, arrows}

\tikzstyle{startstop} = [rectangle, rounded corners, minimum width=3cm, minimum height=1cm,text centered, draw=black, fill=red!30]
\tikzstyle{process} = [rectangle, minimum width=3cm, minimum height=1cm, text centered, draw=black, fill=orange!30]
\tikzstyle{decision} = [diamond, minimum width=3cm, minimum height=1cm, text centered, draw=black, fill=green!30]
\tikzstyle{arrow} = [thick,->,>=stealth]

\usepackage[noindentafter]{titlesec}
\usepackage{bold-extra}

\titlespacing\section{0pt}{14pt plus 4pt minus 2pt}{6pt plus 2pt minus 2pt}
\titlespacing\subsection{0pt}{12pt plus 4pt minus 2pt}{4pt plus 2pt minus 2pt}
\titlespacing\subsubsection{0pt}{8pt plus 4pt minus 2pt}{4pt plus 2pt minus 2pt}
\titlespacing\paragraph{0pt}{6pt plus 4pt minus 2pt}{6pt plus 2pt minus 2pt}

\newcommand{\cF}{c_{\scriptstyle\rm{F}}}
\newcommand{\sF}{s_{\scriptstyle\rm{F}}}
\newcommand{\tF}{\theta_{\scriptstyle\rm{F}}}

\usepackage[T1]{fontenc}
\usepackage{tgtermes}

\newcommand{\carp}{CARPA\xspace}

\title{A Composed Alternating Relaxed Projection Algorithm for Feasibility Problem}
\author{
Yuting Shen\thanks{School of Mathematical Sciences, Shanghai Jiao Tong University, 200240 Shanghai China. E-mail: yutingshen98@sjtu.edu.cn} 
\and 
Jingwei Liang\thanks{School of Mathematical Sciences \& Institute of Natural Sciences, Shanghai Jiao Tong University, 200240 Shanghai China. E-mail: jingwei.liang@sjtu.edu.cn}
}
\date{}

\begin{document}

\setlength{\abovedisplayskip}{5pt}
\setlength{\belowdisplayskip}{3pt}

\maketitle

\begin{abstract}
Feasibility problem aims to find a common point of two or more closed (convex) sets whose intersection is nonempty. In the literature, projection based algorithms are widely adopted to solve the problem, such as the method of alternating projection (MAP), and Douglas--Rachford splitting method (DR). The performance of the methods are governed by the geometric properties of the underlying sets. For example, the fixed-point sequence of the Douglas--Rachford splitting method exhibits a spiraling behavior when solving the feasibility problem of two subspaces, leading to a slow convergence speed and slower than MAP. However, when the problem at hand is non-polyhedral, DR can demonstrate significant faster performance. Motivated by the behaviors of the DR method, in this paper we propose a new algorithm for solving convex feasibility problems. 
The method is designed based on DR method by further incorporating a composition of projection and reflection. A non-stationary version of the method is also designed, aiming to achieve faster practical performance. Theoretical guarantees of the proposed schemes are provided and supported by numerical experiments.

\end{abstract}

\begin{keywords}Feasibility problem, Projection,  Relaxation, Friedrichs angle, Linear convergence
\end{keywords}

\section{Introduction}
\label{sec:introduction}

In this paper, we are interested in the following feasibility problem of two sets 
\beq\label{eq:feasi_problem}
\find \enskip {x\in \calH} \enskip \textrm{such~that} \enskip x \in X \cap Y  ,
\eeq
where $\calH$ is a real Hilbert space, and $X,Y \subset \calH$ are closed convex such that
\[
X \cap Y \neq \emptyset .	
\]
Feasibility problems exist in various fields including signal/image processing \cite{combettes1996convex}, compressed sensing \cite{suantai2019new, van2009algorithm}, seismic data processing \cite{gao2013convergence}, etc. 
In the literature, projection (see Definition \ref{def:projection}) based algorithms are widely adopted to solve feasibility problems. Over the years, significant advances have been made in both numerical methods and theoretical understanding. 

The very basic algorithm to solve the feasibility problem is the method of alternating projection (MAP)~\cite{gubin1967method, bauschke1993convergence}. With an initial point which does not necessarily belong to either of the two sets, MAP applied the projection of the two alternatively. The convergence of the method can be obtained when the two sets are closed and convex, even when the intersection is empty~\cite{bauschke1993convergence}. Based on MAP, various modifications of the method can be found in the literature, such as relaxed alternating projection (RAP) and partial relaxed alternating projection (PRAP)~\cite{bauschke2016optimal}, averaged alternating modified reflections (AAMR) ~\cite{aragon2019optimal}, relaxed averaged alternating reflections (RAAR)~\cite{luke2004relaxed}, simultaneous projection (SP)~\cite{reich2017optimal}, and generalized relaxed alternating projection (GRAP)~\cite{falt2017optimal,dao2018linear}. Recently, by further utilizing the geometric properties of the problem, a so-called circumcentered-reflection method (CRM) is designed \cite{arefidamghani2022circumcentered,behling2021circumcentered,behling2022successive}. 
Operator splitting methods can also be applied to solve the feasibility problem, since the projection operator is a special case of the proximity operator \cite{bauschke2011convex}. Actually, MAP is an instance of the Backward-Backward splitting method \cite{Bauschke2005}. For both AAMR and GRAP methods, with proper parameter choices, they can recover the Peaceman--Rachford splitting method \cite{peaceman1955numerical}, and the Douglas--Rachford splitting method \cite{douglas1956numerical}; see Table \ref{tab:rates} for details.

Theoretical understandings of projection based methods under different settings are also studied. As most projection based methods can be written fixed-point iteration of nonexpansive operators, their global convergence and convergence rate can be found in \cite{liang2016convergence}. 
When the underling two sets are so-called ``partly smooth'' \cite{Hare2007}, locally linear convergence can be obtained, see for instance \cite{Lewis2009,liang2017localDR}. In \cite{Li2016}, convergence properties of DR for nonconvex feasibility problem is studied. 
When the two sets are moreover subspaces, in \cite{bauschke2014rate,Bauschke2016,aragon2019optimal} the authors studied the (optimal) linear convergence rate of several projection based methods, including DR and MAP. We also refer to \cite{Dao2019} for linear convergence of projection methods, and \cite{Bauschke2023,Bauschke2024} for recent development of projection algorithms.

\subsection{Motivation and contribution}

\begin{figure}[!htb]
	\centering
 \includegraphics[width=0.5\linewidth, trim={0mm 4mm 0mm 2mm}, clip]{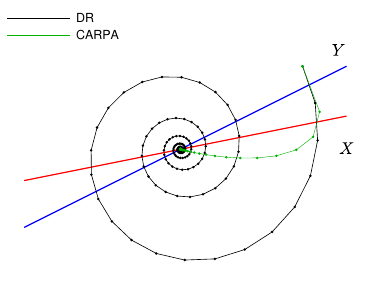}
	\caption{Spiraling behavior of the fixed-point sequence of Douglas--Rachford splitting method.}
	\label{fig:subbehavior}
\end{figure}

Consider the problem where $X, Y$ are two subspaces, then MAP and DR converge linearly where the rate of convergence is characterized by the Friedrichs angle between $X$ and $Y$ (see Definition \ref{def:The Friedrichs angle}) \cite{bauschke2014rate,bauschke2016optimal}. 
Moreover, from the convergence analysis of \cite{bauschke2014rate}, it can be derived that the trajectory of the fixed-point sequence of DR forms a spiral, see also \cite{poon2019trajectory}. 
In Figure \ref{fig:subbehavior} we provide an illustration of this behavior, $X, Y$ are two intersecting subspaces, the {\it black spiral} denotes the trajectory of the fixed-point sequence of DR. 
Such spiral behavior is the main reason why DR performs poorly on subspace feasibility problems, see Table \ref{tab:rates} and Figure \ref{fig:W1comparenew} for details. 

%%%%%%%%%%%%%%%%%%%%%%%%%%%%%%%%%%%%%%%%%%%%%%%%%%%%%%%%

In \cite{bauschke2016optimal}, it is shown that for this problem the standard DR achieves the optimal convergence rate, and relaxation does not provide acceleration. Moreover, as pointed out in \cite{poon2019trajectory}, inertial acceleration also fails to provide speedup. Motivated by this behavior of DR method, in this paper we propose a composed alternating relaxed projection algorithm (\carp, see Algorithm \ref{alg:cdr}) to improve the performance of the DR method. The {\it green} trajectory in Figure \ref{fig:subbehavior} displays the trajectory of the fixed-point sequence of \carp. The key idea of the method a convex combination of the DR and projection-reflection, doing so the method managed to provide significant acceleration over DR method, see Table \ref{tab:rates} and Figure \ref{fig:W1comparenew}. More precisely, this work contains the following contributions
\begin{itemize}
    \item A novel projection based algorithm (Algorithm \ref{alg:cdr}) for solving convex feasibility problems, with theoretical guarantees. We also design a non-stationary version of the method, and prove its convergence property.

    \item For the feasibility problem of two subspaces, we analyze the linear rate of convergence of the method and identify parameter choices under which optimal convergence rate is obtained. 

    % \item 
\end{itemize}
Furthermore, numerical experiments on feasibility of two subspaces, a ball tangent to a subspace and a compressed sensing problem are provided to demonstrate the advantages of the proposed schemes.

\subsection{Paper organization} 
The remainder of this paper is organized as follows. In Section \ref{sec:maths_bg}, we introduce the notations and basic results used throughout the paper. Section \ref{sec:carpa} describes the proposed method and its non-stationary version, with theoretical analysis. 
In Section \ref{sec:numerics}, we provide numerical comparisons against several existing methods in the literature.

\section{Mathematical background}
\label{sec:maths_bg}

Denote $\calH$ the Hilbert space equipped with scalar inner product $\iprod{\cdot}{\cdot}$ and the associated norm $\norm{\cdot}$.  
Denote $\bbR^{n}$ the $n$-dimensional real Euclidean space.
{$\Id$ denotes the identity operator.} 
Given $A\in\bbR^{m\times n}$, %$A^*$ stands for the adjoint (complex transposed) matrix of A. We 
denote by $ \Ker (A)$, $\ran (A)$, $\rank (A)$ the kernel, range, rank of $A$, respectively. % The set $\Fix (A):= \ker(A -\Id)$ is known as the set of fixed-points of A. 

\paragraph{Operators} In this part, we provide the definitions of nonexpansive and projection related operator, we refer to \cite[Chapter 4]{bauschke2011convex} for more details regarding this topic. 

\begin{definition}[Nonexpansive operator]\label{dfn:nonexpan-operator}
An operator $\calF : \calH\rarrow\calH$ is called {nonexpansive} if it is $1$-Lipschitz continuous, \ie
\beqn
\norm{\calF  (x) - \calF  (x')} \leq \norm{x-x'} ,~~~ \forall x,x' \in\calH .
\eeqn
For any $\alpha \in ]0,1[$, $\calF $ is $\alpha$-averaged if there exists a nonexpansive operator $\calF '$ such that
\[
\calF = \alpha \calF '+(1-\alpha)\Id .
\]
When $\alpha=1/2$, $\calF$ is also called ``firmly nonexpansive'', and $2\calF-\Id$ is nonexpansive in this case. 
\end{definition}

A point $x$ is called a fixed-point of $\calF$ if it satisfies $x = \calF(x)$, and the set of fixed-points of $\calF$ is defined as $\fix(\calF) = \Ba{ x\in\calH \mid x = \calF(x) }$, which could be empty. 

\begin{definition}[Projection]\label{def:projection}
Let set $X \subset \calH$ be nonempty closed convex, its {projection operator} is defined by, %$\forall x \in \calH$
\[
\proj_{X}(x)
\eqdef \Argmin_{y \in X}~ \norm{x-y} ,~~ \forall x \in \calH ,
\]
which is firmly nonexpansive. 

\end{definition}

\begin{definition}[Relaxed projection and reflection]\label{def:relax_projection}
Let set $X \subset \calH$ be nonempty closed convex and $r \in [-1, 1]$, the relaxed projection is defined by, %$\forall x \in \calH$
\[
\scrR^{r}_{X}(x)
\eqdef \Pa{ (1+r) \proj_{X} - 
r\Id} (x), ~~~ \forall x \in \calH
\]
which is $\frac{1+\gamma}{2}$-averaged. 
When $r=1$, $\scrR^{r}_{X}$ is called the {``reflection operator''} and is denoted by $\scrR_{X}$ for simplicity. 

\end{definition}

\paragraph{Angles between subspaces} 
Let $X, Y \subset \bbR^n$ be two subspaces, we have the following definitions regarding the angles between them \cite{bauschke2014rate}.

\begin{definition}[Principal angles] \label{def:principal angle}
% \cite{bauschke2016optimal,meyer2023matrix}
The principal angles $\theta_k \in$ $\left[0, \frac{\pi}{2}\right], k=1, \cdots, p$ between two subspaces $X, Y \subset \bbR^n$ where $p = {\min}\Ba{ \dim (X), \dim (Y)}$ are recursively defined by
$$
\begin{aligned}
\cos (\theta_k) 
& =\left\langle u_k, v_k\right\rangle \\
& =\max \left\{\begin{array}{l}
\langle u, v\rangle \left\lvert\, 
\begin{array}{l}
u \in X, v \in Y,\|u\|=\|v\|=1, \\
\left\langle u, u_j\right\rangle=\left\langle v, v_j\right\rangle=0, j=1, \ldots, k-1
\end{array}\right.
\end{array}\right\} \qwithq u_0=v_0:=0 .
\end{aligned}
$$
It is worth mentioning that the vectors $u_k, v_k$ are not uniquely defined, but the principal angles $\theta_k$ are unique with $0 \leq \theta_1 \leq \theta_2 \leq \cdots \leq \theta_p \leq \frac{\pi}{2}$. 
\end{definition}

\begin{definition}[Friedrichs angle]\label{def:The Friedrichs angle}
  The cosine of the Friedrichs angle  between $X$ and $Y$ is
$$
\cF (X, Y)= \cos (\tF) :=\sup \left\{|\langle u, v\rangle| \mid u \in X \cap(X \cap Y)^{\bot}, v \in Y \cap(X \cap Y)^{\bot},\|u\| \leq 1,\|v\| \leq 1\right\} .
$$
\end{definition}
% We write $\cF (X, Y)$ for $\cF $ if we emphasize the subspaces utilized.
We will use $\cF$ for simplicity if there is no confusion.

\begin{proposition}[Principal angles and Friedrichs angle]
    Let $s := \dim(X \cap Y)$, then $\theta_k = 0 \  {\rm for}\ k=1,\ldots,s $ and $\theta_{s+1}=\tF>0.$
\end{proposition}

\paragraph{Convergence rate of matrices}

For any $A \in \bbR^{n\times n}$, $\sigma(A)$ and $\rho(A)$ denote the spectrum and spectral radius of $A$, respectively. 
% which means the set of all eigenvalues of $A$.The spectral radius of A is denoted by $\rho(A)$. 
We say $A$ converges linearly to $A^{\infty}$ if there exists $r\in [0,1[$ and $M,N > 0$ such that
\beqn
\|A^k - A^\infty \| \leq M r^k,~~~ \forall k \geq N . 
\eeqn
An eigenvalue $\lambda \in \sigma(A)$ is {\it semisimple} if 
$$
\ker (A-\lambda \Id) = \ker (A -\lambda \Id)^2.
$$  
The {\it subdominant eigenvalue} $A$ is defined by
$$
\xi(A):=\max \{|\lambda| \mid \lambda \in\{0\} \cup \sigma(A) \backslash\{1\}\}, 
$$
\ie the largest eigenvalue besides $1$ without sign.

\begin{fact}[{\cite[Fact 2.4]{bauschke2016optimal}}]
Suppose that $ A \in \bbR^{n \times n}$, $\lim_{k \rightarrow \infty} A^k$ exists if and only if
\beqn
\rho(A) < 1 \qorq \rho(A) = 1 \qandq\lambda = 1 \text{ is semisimple and it is the only eigenvalue on the unit circle}.
\eeqn
\end{fact}

\begin{theorem}[{\cite[Theorem 2.12, Theorem 2.15]{bauschke2016optimal}}]
Suppose that $ A \in \bbR^{n \times n}$ is convergent to $A^{\infty} \in \bbR^{n \times n}$, then  $\xi(A)=\rho\left(A-A^{\infty}\right)<1$ and that
$$
\left(A-A^{\infty}\right)^k=A^k-A^{\infty} \quad \text { for all } k \in \mathbb{N} .
$$ 
Moreover, the followings hold 
\begin{itemize}
    \item[(i)] $A$ is linearly convergent with any rate $ r \in ]\xi(A), 1[$.
    \item[(ii)] $\xi(A)$ is the optimal convergence rate of A if and only if all the subdominant eigenvalues are semisimple.
\end{itemize}
\end{theorem}

\section{A composed alternating relaxed projection algorithm}\label{sec:carpa}

This section is dedicated to the formulation and analysis of our proposed algorithm. We first provide the details of the method with convergence analysis, and then discuss a non-stationary version of the method, which demonstrates more practical performances. 

\subsection{The proposed algorithm}
Motivated by the drawback of DR method for solving (locally) polyhedral feasibility problems, we design a new algorithm called the ``composed alternating relaxed projection algorithm (\carp)'', whose details are provided below.

\begin{center}
	\begin{minipage}{0.95\linewidth}
		\begin{algorithm}[H]
			\caption{A composed alternating relaxed projection algorithm (\carp)} \label{alg:cdr}
			{\noindent{\bf{Initial}}}: $k = 0$, $z^{(0)} \in \calH, x^{(0)} = \proj_{X} (z^{(0)})$, $\gamma \in [0, 1[$, $\mu \in ]0, 2/(1+\gamma)[$\;
			\Repeat{convergence}{
				\beq\label{eq:cdr}
				\begin{aligned}
					\xkp &= \proj_{X} (\zk) , \\
					\ukp &= 2\xkp - \zk , \\
					\ykp &= \proj_{Y}\pa{\ukp} , \\
					\zkp &=  (1-\mu)\zk + \mu\Pa{  (1-\gamma)( \zk +\ykp - \xkp) + \gamma \ykp } . \\
				\end{aligned}
				\eeq
				$k = k + 1$\;
			}
		\end{algorithm}
	\end{minipage}
\end{center}

In what follows, we discuss the convergence properties of \carp, starting from analyzing the fixed-point formulation of the method and then convergence of the sequence. 
Define the following operators
\beq\label{eq:fcarp}
\calF_{1} = \sfrac{1}{2}\pa{\Id + \scrR_{\calY}\scrR_{\calX}}, \quad \calF_{2} = \proj_{\calY}\scrR_{\calX} . % \\
\eeq
Then we denote 
\[
\fCDR 
= (1-\gamma) \calF_{1} + \gamma \calF_{2}
\qandq
\fCDRmu 
= (1-\mu)\Id + \mu \fCDR . 
\]
We have the following proposition on the properties of these operators. 

\begin{proposition}
The operator $\fCDR$ is $\frac{1+\gamma}{2}$-averaged nonexpansive, and $\fCDRmu$ is $\frac{(1+\gamma)\mu}{2}$-average. 
Moreover, the proposed algorithm can be written as the fixed-point iteration of $\zk$ which reads
\beq\label{eq:fp_cdr}
\zkp = \fCDRmu (\zk) .
\eeq 
\end{proposition}
\begin{proof}
From the iteration of \carp, we have
\[ 
\ykp 
= \proj_Y \pa{2\proj_X - \Id} (\zk) 
= \calF_{2} (\zk) 
\]
and 
\[
\begin{aligned}
\zk +\ykp - \xkp
&= \Pa{\Id + \proj_Y \pa{2\proj_X - \Id} - \proj_X} (\zk) \\
&= \sfrac{1}{2}\Pa{ \Id + \scrR_Y \scrR_X } (\zk) \\
&= \calF_1(\zk) .
\end{aligned}
\]
Combine the two formulas we get $\zkp = \fCDRmu(\zk)$. 

Back to the expression of $\fCDR$, we have
\[
\begin{aligned}
 \fCDR 
 &= (1-\gamma) \calF_{1} + \gamma \calF_{2}  \\
 &= (1-\gamma) \sfrac{1}{2}\pa{\Id + \scrR_{\calY}\scrR_{\calX}} + \gamma \proj_{\calY}\scrR_{\calX}  \\
 % &= \sfrac{1-\gamma}{2} \Id + \sfrac{1}{2} \scrR_{\calY}\scrR_{\calX} + \gamma \Pa{ \proj_{\calY}\scrR_{\calX} - \sfrac{1}{2} \scrR_{\calY}\scrR_{\calX} }    \\
 % &= \sfrac{1-\gamma}{2} \Id + \sfrac{\gamma}{2}\scrR_{\calX} +   \sfrac{1}{2}\scrR_{\calY}\scrR_{\calX}  \\
 &= \sfrac{1-\gamma}{2} \Id +  \sfrac{1+\gamma}{2} \Pa{ \sfrac{1-\gamma}{1+\gamma}\scrR_{\calY} + \sfrac{2\gamma}{1+\gamma} \proj_{Y} } \scrR_{\calX}  .
\end{aligned}
\]
Owing to the result of \cite[Chapter 4]{bauschke2011convex}, $\proj_Y$ is firmly nonexpansive and $\scrR_Y$ is nonexpansive, hence $\sfrac{1-\gamma}{1+\gamma}\scrR_{\calY} + \sfrac{2\gamma}{1+\gamma} \proj_{Y}$ is nonexpansive \cite[Proposition 4.6]{bauschke2011convex}. Consequently, $\Pa{ \sfrac{1-\gamma}{1+\gamma}\scrR_{\calY} + \sfrac{2\gamma}{1+\gamma} \proj_{Y} } \scrR_{\calX}$ is nonexpansive. By the definition of averaged operators, we conclude the proof. 
The averagedness of $\fCDRlam$ follows immediately. %is $\frac{(1+\gamma)\mu}{2}$-averaged. % with $\lambda = 1 - \mu$. % (see Definition \ref{dfn:nonexpan-operator}). 
\end{proof}

\begin{proposition}[Set of fixed-points]\label{lem:fxs_cdr}
Let $X,Y\subset\calH$ be nonempty closed convex sets with nonempty intersection, \ie $X\cap Y \neq \emptyset$, then $\fix(\fCDR) = X \cap Y$. 
\end{proposition}

% This result is quite  the sake of 

\begin{proof}
We provide a simple proof for the sake of completeness. 
Owing to \cite[Proposition 4.47]{bauschke2011convex}, we have
\[
\fix(\fCDR) = \fix(\calF_{1}) \cap \fix(\calF_{2}) .
\]
We first characterize the set of fixed-points for $\calF_{2}$, for which we have
\[
\fix(\calF_2) = X \cap Y .
\]
It is obvious that for any $z\in X\cap Y$, $z = \calF_2(z)$. 
Suppose there exists $z\in\calF_2(z)$ such that $z\notin X\cap Y$. Now that $z\in\calF_2(z)$, it means
\[
z = \proj_Y \scrR_X (z) \in Y ,
\]
which further implies
\[
\begin{aligned}
\scrR_X - z \in N_{Y}(z)
&\quad\Longleftrightarrow\quad
2\proj_X(z) - 2z \in N_{Y}(z)  \\
&\quad\Longleftrightarrow\quad
\proj_X(z) - z \in N_{Y}(z)  \\
&\quad\Longleftrightarrow\quad
z = \proj_Y \proj_X (z) ,
\end{aligned}
\]
where $N_{Y}(\cdot)$ denotes the normal cone operator. 
% As $z\notin X\cap Y$, this means $z \notin X$, and  
% \[
% \proj_X(z) \in X \setminus (X \cap Y)  . 
% \]
As a result, we have $z\in X\cap Y$ \cite{bregman1965method} which conflicts with the assumption, hence we obtain $\fix(\calF_2) = X \cap Y$. Now for any $z \in X \cap Y$, it is easy to verify that 
\[
z \in \fix(\calF_1) ,
\]
and we conclude the proof. \qedhere
\end{proof}

With above characterization of the set of fixed-points, we are ready to provide the convergence result of our proposed algorithm.

\begin{proposition}[Convergence of \carp]\label{thm:convergence_cdr}
	For problem \eqref{eq:feasi_problem} and Algorithm \ref{alg:cdr}, suppose $X\cap Y\neq\emptyset$ and let $\gamma \in ]0, 1[$ and $\mu\in]0,2/(1+\gamma)[$, then there exists a $\zsol \in \fix(\fCDR)$ such that $(\xk,\yk,\zk) \to (\zsol,\zsol,\zsol)$. 
\end{proposition}
\begin{proof}
As our proposed algorithm can be written as a fixed-point iteration in terms of $\zk$(see \eqref{eq:fp_cdr}), we have the convergence of $\zk$ due to  \cite[Theorem 5.15]{bauschke2011convex}, that is 
\[
\zk \to \zsol \in X \cap Y .
\]
Given this, we have 
\[
\xk \to \xsol = \proj_{X}(\zsol) = \zsol \qandq
\yk \to \ysol = \proj_{Y}(2\xsol-\zsol) = \proj_{Y}(\zsol) = \zsol ,
\]
which concludes the proof. \qedhere
\end{proof}

\begin{remark}[Relation with existing work]
In this part, we briefly discuss the relation between our proposed method and several existing ones in the literature. It is obvious that when $\gamma=0$, our method recovers the relaxed DR scheme, and $\gamma=1$ the relaxed projection-reflection method which can be seen as a special case of \cite{falt2017optimal,dao2018linear,aragon2019optimal}. 
When replacing the $\calF_2$ operator in \eqref{eq:fcarp} with $\proj_X$, one recovers the method proposed in \cite{luke2004relaxed}. 
\end{remark}

\subsection{A non-stationary version of \carp (ns\carp)}

In this section, we consider a non-stationary version of Algorithm \ref{alg:cdr} by fixing $\mu$ and changing the value of $\gamma$ over iterations, aiming to achieve practical acceleration. The proposed scheme is described below in Algorithm~\ref{alg:adcdr_adaptive1}. 

\begin{center}
	\begin{minipage}{0.95\linewidth}
		\begin{algorithm}[H]
			\caption{A non-stationary version of \carp}  \label{alg:adcdr_adaptive1}
			{\noindent{\bf{Initial}}}: $k = 0$, $z^{(0)} \in \bbR^n, x^{(0)} = \proj_{X} (z^{(0)})$, $\gamma_0,\gamma_{\max}, \gamma_{\min} \in [0, 1[$, $\mu \in ]0, 1]$, and $c_1 ,c_2 ,\delta >0 $\;
			\Repeat{convergence}{
				Compute the updates:
				\begin{equation*}
				\begin{aligned}
					\xkp &= \proj_{X} (\zk) , \\
					\ukp &= 2\xkp - \zk , \\
					\ykp &= \proj_{Y}\pa{\ukp} , \\
					\zkp &= (1-\mu)\zk + \mu\Pa{ \zk + (1-\gmk)(\ykp - \xkp) + \gmk(\ykp - \zk) } , \\[1mm]
				\end{aligned}
				\end{equation*}
                
				If \(k > 0\), let \( \rho_k = \frac{\|\zkp - \zk\|}{\|\zk - \zkm \|} \), then \vskip2pt
                
                \beq\label{eq:gammak}
                \begin{aligned}
                \gamma_{k+1/2} &= \left\{ \begin{aligned} 
                    \gmk + \sfrac{c_2}{(k+1)^{2+\delta}}: & \quad \rho_k < c_1 , \\
                    \gmk - \sfrac{c_2}{(k+1)^{2+\delta}}: & \quad o.w. 
                    \end{aligned} \right.  \\[2mm]
                    \gmkp &= \max\bBa{\min\Ba{\gamma_{k+1/2},~ \gamma_{\max}},~\gamma_{\min}} . 
                    \end{aligned}
                    \eeq
                $k = k + 1$\;
			}
		\end{algorithm}
	\end{minipage}
\end{center}

\begin{remark}
\label{alg:adcdr_adaptive2}
One can also consider updating $\gmk$ via the Armijo--Goldstein-like rule:
\beqn
\gmk =\left\{
\begin{array}{ll}
  \gmk \cdot \eta, &\rho_k < c_1 \\ 
  \min\{ \gmk / \eta, 1 \}, & o.w.  
\end{array} 
\right.
\eeqn
% $\seq{\gmk}$ finally converges to $\gmsol$ or enters to a periodic state. 
% 
% When within the periodic state, the convergence of ns\carp is also can be proved easily. 
% Besides, 
Other choices include: let $\gamma\in]0,1[$
\[
\gmk = \gamma + \sfrac{1}{k^\delta},~\delta>1    ,\quad
\gmk = \gamma + \rho^k,~\rho\in ]0,1[      \qandq 
\gmk = \gamma + \exp(-{k}/{\delta}),~\delta > 0 . 
\]
Note that for the above three choices, we have $\gamma_k\to\gamma$ and $\seq{\gamma_k-\gamma}$ is summable. 
\end{remark}

\begin{remark}\label{rmk:nsdr}
A non-stationary version of DR (nsDR) proposed in \cite{lorenz2019non} will be included in the following comparisons. 
    \beq\label{eq:nsdr}
    \begin{aligned}        
					\xkp &= \proj_{X} (\zk) , \\
                        \tau_k & = \sfrac{\norm{\xkp}}{\norm{\xkp - \zk}} \\
                        \ukp & = (1+\tau_k)\xkp-\tau_k \zk\\
					\ykp &= \proj_{Y}\pa{\ukp} , \\
					\zkp &= \ykp + \tau_k(\zk - \xkp). \\
    \end{aligned}
    \eeq
\end{remark}

Following the result of \cite{liang2016convergence}, we can analyze the convergence of the above scheme by treating the non-stationarity as a perturbation error. More precisely, Algorithm \ref{alg:adcdr_adaptive1} can be written as
\begin{equation}
\label{eq:fp_nscdr}
\zkp 
= \calF_{\gmk, \mu}(\zk)
= \calF_{\gamma, \mu}(\zk) + \pa{ \calF_{\gmk, \mu} - \calF_{\gamma, \mu} }(\zk) ,
\end{equation}
where $\pa{ \calF_{\gmk, \mu} - \calF_{\gamma, \mu} }(\zk)$ is the perturbation error. 
Note that this error is controlled by the updating rule of $\gamma_k$, hence properly designed rules can ensure the convergence of the scheme, this is the reason for choosing \eqref{eq:gammak} and the content of the proposition below.

\begin{proposition}[Convergence of ns\carp]
For problem \eqref{eq:feasi_problem} and Algorithm \ref{alg:adcdr_adaptive1}, suppose $X\cap Y\neq\emptyset$. 
With the updating rule \eqref{eq:gammak}, there exists a $\zsol \in \fix(\fCDR)$ such that $(\xk,\yk,\zk) \to (\zsol,\zsol,\zsol)$. 
\end{proposition}
 
The situation here is simpler than that of \cite{liang2016convergence}, as we only have $\gamma_k$ varying. If we can show that $\pa{ \calF_{\gmk, \mu} - \calF_{\gamma, \mu} }(\zk)$ is summable, then by virtue of quasi-Fej\'er monotonicity \cite{bauschke2011convex}.  

\begin{proof}
Note that if the perturbation error
\[
\ek = \pa{ \calF_{\gmk, \mu} - \calF_{\gamma, \mu} }(\zk)
\]
satisfies 
\[
\sum_{k=1}^{\infty} \norm{\ek} < \infty . 
\]
Then the sequence $\seq{\zk}$ is quasi-Fej\'er monotone and its convergence follows easily \cite{bauschke2011convex}.  To show this, we have the following steps
\begin{itemize}

\item[i).] We first study the properties of $\seq{\gamma_k}$. By the definition of $\gamma_k$, we have $|\gamma_{k+1} - \gamma_k| \leq \frac{c_2}{(k+1)^{2+\delta}}$.
Given any $\epsilon>0$, let $N = \lceil\big(\frac{c_2}{(1+\delta)\epsilon} \big)^{\frac{1}{1+\delta}} \rceil $, then $j > i > N$, we have
\beq
\label{mn}
\begin{aligned}
|\gamma_{j} - \gamma_{i}| 
% &= | \gamma_j - \gamma_{j-1} + \cdots \gamma_{i+1} - \gamma_i| \\
\leq | \gamma_j - \gamma_{j-1} | +\cdots +|\gamma_{i+1} - \gamma_i| 
&=\sum_{\ell=i}^{j-1} \frac{c_2}{(\ell+1)^{2+\delta}}  \\ 
&\leq \int_i^\infty \frac{c_2}{x^{2+\delta}}dx \\
&= \frac{c_2}{(1+\delta)i^{1+\delta}} \\
&< \frac{c_2}{(1+\delta)N^{1+\delta}} 
< \epsilon .
\end{aligned}
\eeq
% \rc{We don't need the above property...?}
This means $\seq{\gamma_k}$ is a Cauchy sequence, hence convergent.
Denote $\gmsol$ the limit of $\gamma_k$, and from \eqref{mn} we have 
\beqn
|\gamma_k-\gmsol|
\leq \sum_{i=k}^{\infty} \frac{c_2}{(i+1)^{2+\delta}}
\leq \frac{c_2}{(1+\delta)k^{1+\delta}} .
\eeqn
Consequently,
\beqn
\begin{aligned}
\sum_{k=1}^{\infty}\mu|\gamma_k -\gmsol| 
\leq \sum_{k=1}^{\infty}\frac{\mu c_2}{(1+\delta)k^{1+\delta}} 
&=  \frac{\mu c_2}{1+\delta} \sum_{k=1}^{\infty}\frac{1}{ k^{1+\delta}} 
< \infty  .
\end{aligned}
\eeqn

\item[ii).] Next we show that the sequence $\seq{\zk}$ is bounded. 
As for all $k$, we have $\gamma_k \in [0,1[$, $\calF_{\gmk, \mu}$ is nonexpansive and $X\cap Y \subset \fix(\calF_{\gmk, \mu})$. Let $\zsol \in X \cap Y$, then 
\[
\begin{aligned}
\norm{\zkp - \zsol}
= \norm{ \calF_{\gmk, \mu} (\zk - \zsol) }
&\leq \norm{\zk - \zsol } \\
&\leq \norm{z^{(0)}-\zsol } 
\end{aligned}
\] 
which means $\seq{\zk}$ is bounded. Denote
\[
r = \sup_{k\in\bbN} \norm{\zk} .
\]
Let $\gamma=\gamma^\star$, by the definition of $\calF_{\gmk, \mu}$ we have
\[
\begin{aligned}
\calF_{\gmk, \mu} - \calF_{\gamma, \mu}
&= \Pa{(1-\mu)\Id + \mu (1-\gamma_k) \calF_{1} + \mu \gamma_k \calF_{2}} - \Pa{(1-\mu)\Id + \mu (1-\gamma) \calF_{1} + \mu \gamma \calF_{2}}  \\
&= -\mu(\gamma_k-\gamma) \calF_{1} + \mu (\gamma_k-\gamma)\calF_2 .
\end{aligned}
\]
Recall that both $\calF_1,\calF_2$ are nonexpansive, then from \eqref{eq:fp_nscdr} we get: 
\[
\begin{aligned}
\norm{\zkp-\zsol}
% = \calF_{\gmk, \mu}(\zk)
&\leq \norm{ \calF_{\gamma, \mu}(\zk) - \zsol} + \norm{ \pa{ \calF_{\gmk, \mu} - \calF_{\gamma, \mu} }(\zk - \zsol)}   \\
&\leq \norm{ \calF_{\gamma, \mu}(\zk) - \calF_{\gamma, \mu}(\zsol)} + 2\mu \abs{\gamma_k-\gamma} \norm{\zk-\zsol} \\
&\leq \norm{ \zk - \zsol} + 4r\mu \abs{\gamma_k-\gamma}  .
\end{aligned}
\]
Since $\abs{\gamma_k-\gamma}$ is summable, we get that $\seq{\zk}$ is quasi-Fej\'er monotone. 

\end{itemize}
By virtue of \cite[Proposition 5.34]{bauschke2011convex} we obtain the convergence of the sequence $\seq{\zk}$, and consequently the sequences $\seq{\xk,\yk}$. \qedhere
\end{proof}

\section{Numerical experiments}\label{sec:numerics}

In this section, we demonstrate the performance of the proposed methods by considering several feasibility problems, including two subspaces, a line tangent to a ball, and compressed sensing problems. 

\subsection{Feasibility of two subspaces}\label{rate}

We first consider a feasibility problem of two subspaces $X,Y\subset \bbR^n$, with $p = \dim(X)$ and $q =\dim(Y)$.
Denote $\tF$ the Friedrichs angle between $X$ and $Y$, and $\theta_p$ the largest principal angle. Furthermore, denote $X^\bot, Y^\bot$ the orthogonal complement of $X$ and $Y$, respectively.

\paragraph{Optimal convergence rates}\label{otherprojection_methods}

As there is a rich literature of projection based methods for solving feasibility problems, in Table \ref{tab:rates} we summarize several representative ones and their optimal convergence rate for the considered problem.

\begin{table}[!htb]
\centering  
\begin{tabular}{|m{6.2cm}|m{6cm}|m{2.8cm}|}
\hline
Method & Fixed-point operator  & Optimal rate  \\ \hline
Simultaneous projection (SP) \cite{reich2017optimal} & $\frac{1}{2}(\proj_Y + \proj_X)$ & $\frac{1+ \cos(\tF)}{2}$ \\ \hline
Relaxed alternating projection (RAP) \cite{bauschke2016optimal} & $(1-\mu) \Id + \mu \proj_Y \proj_X$ & $\sfrac{1-\sin^2(\tF)}{1+\sin^2(\tF)}$ \\ \hline
Partial relaxed alternating projection (PRAP) \cite{bauschke2016optimal} & $(1-\mu)\proj_Y +\mu \proj_Y \proj_X$ & $\frac{\sin^2(\theta_p)-\sin^2(\tF)}{\sin^2(\theta_p)+\sin^2(\tF)}$ \\ \hline
Generalized relaxed alternating projection (GRAP) \cite{falt2017optimal,dao2018linear} & $(1-\mu)\Id + \mu \scrR^{\alpha_2}_Y \scrR^{\alpha_1}_X$ & $\frac{1-\sin(\tF)}{1+\sin(\tF)}$ \\ \hline
Averaged alternating modified reflections (AAMR) \cite{aragon2019optimal} & $(1-\mu) \Id + \mu(2\beta \proj_Y - \Id)(2\beta \proj_X - \Id )$ & $\frac{1-\sin(\tF)}{1+\sin(\tF)}$ \\ \hline
Relaxed averaged alternating reflections (RAAR) \cite{luke2004relaxed} & $\mu(\proj_Y \proj_X+\proj_{Y^{\bot}} \proj_{X^{\bot}}) + (1-\mu)\proj_X$ & $\frac{\cos (\tF)}{\ssqrt{1+2\cos (\tF) \sin (\tF)}}$ \\ \hline
A relaxed Douglas--Rachford (DRAP) \cite{thao2018convergent} & $\proj_Y \proj_X + \mu \proj_{Y^{\bot}} \proj_{X^{\bot}}$ & $1-\sin (\tF)$ \\ \hline
% 
 % &  &  \\ \hline
 % &  &  \\ \hline
\end{tabular}
\caption{Optimal convergence rates for two-subspace feasibility problem for $p=q=n/2$. The corresponding parameter choices are: 1) RAP $\mu=\frac{2}{1+\sin^2(\tF)}$; 2) PRAP $\mu = \tfrac{2}{\sin^2(\theta_p)+\sin^2(\tF)}$; 3) GRAP $\mu =1$ and $\alpha_1 = \alpha_2 = \tfrac{1-\sin(\tF)}{1+\sin (\tF)}$; 4) AAMR $\mu = 1$ and $\beta = \tfrac{1}{1+\sin(\tF)}$; 5) RAAR $\mu=\tfrac{ 1}{1+2\cos (\tF) \sin (\tF) }$; 6) DRAP $\mu=\tfrac{(1-\sin( \tF))^2}{\cos^2(\tF)}$. }
\label{tab:rates}
\end{table}

\begin{remark}%$~$
Recall that for two subspaces feasibility problem, the convergence of MAP is $\cos^2(\tF)$, while for DR it is $\cos(\tF)$ \cite{bauschke2016optimal}.
Note that GRAP and AAMR have the same optimal rate, as they use different representations for relaxed projection operator. 
For RAAR and DRAP, their optimal convergence rates are not previously available to the literature, and we refer to Appendix \ref{drapraar} for brief derivations. 
\end{remark}

\begin{figure}[!htb]
	\centering
\includegraphics[width=0.6\linewidth]{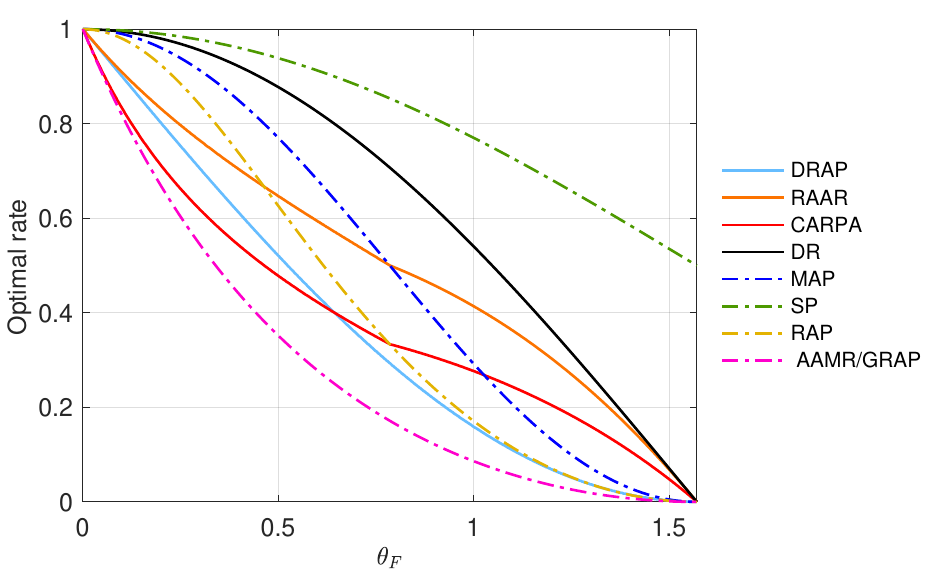}\\[-2mm]
	\caption{Optimal convergence rates of several representative methods.}
	\label{fig:W1comparenew}
\end{figure}

Since all methods can be written as fixed-point iteration, our comparison criterion is $\norm{\zk-\zsol}$, where $\zk$ represents the fixed-point sequence and $\zsol$ denotes a fixed-point of each method. 
In Figure \ref{fig:W1comparenew} we provide a graphical illustration of the optimal convergence rates of our proposed scheme \carp, MAP/DR and the methods in Table \ref{tab:rates} for $\tF\in[0,\pi/2]$, $\theta_p=\pi/2$ and $p=q$. 
We refer to Appendix \ref{pf:Fgammamu} for a detailed derivation of the optimal convergence rate of the proposed \carp scheme. 
% In Figure \ref{fig:W1comparenew}, we provide a comparison of the optimal convergence rate of the aforementioned methods when $\theta_p=\pi/2$.
% \todo{$p=q, p+q=n$ compute by numerical experiments (RAAR, DRAP, CARPA all with relaxeion $\mu$ )} 
% Note that our method admits at least second best optimal convergence rate. 
From the comparison, we observe that
\begin{itemize}
    \item AAMR and GRAP have the best optimal rate over all methods.
    \item For our proposed \carp, it is the second best when $\tF$ is relatively small (smaller than about $\pi/5$). While for $\tF\in [\pi/5, \pi/2]$, \carp is slower than AAMR/GRAP, DRAP, RAP and MAP. 
    % \rc{AP or MAP?}
\end{itemize}

\paragraph{Numerical illustrations}

To demonstrate the above results, below we provide a simple experiment. %for different purposes. 
The basic settings are: $n=100$, $\dim(X)=\dim(Y)=50$ such that $X\cap Y = \Ba{0}$. Moreover, we fix $\theta_p = \pi/2$ and four choices of $\tF$,
\[
\tF \in \Ba{ 0.1,~ 0.4,~ 0.7,~ 1 } .
\]
For all methods, we compare the absolute error of fixed-point sequence to a fixed-point of the corresponding method. The stopping tolerance is $10^{-12}$.

\begin{figure}[!htb]
    \centering
    \subfloat[{$\tF=0.1$}]{\includegraphics[width=0.425\textwidth]{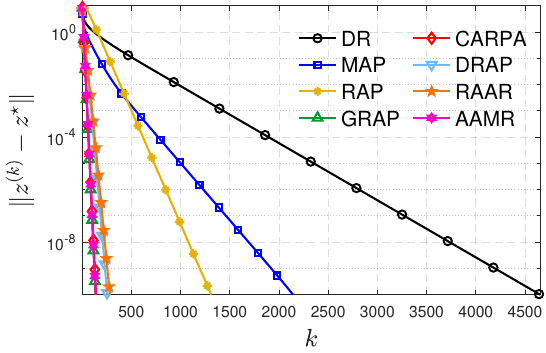}\label{fig:Experiment2-1}}
    \hskip3mm
    \subfloat[{$\tF=0.4$}]{\includegraphics[width=0.425\textwidth]{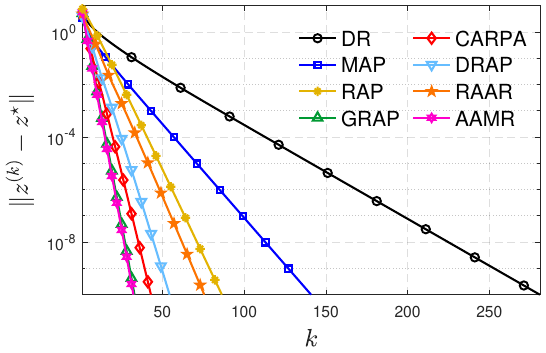}\label{fig:Experiment2-2}}
    \\[-3mm]
    \subfloat[{$\tF=0.7$}]{\includegraphics[width=0.425\textwidth]{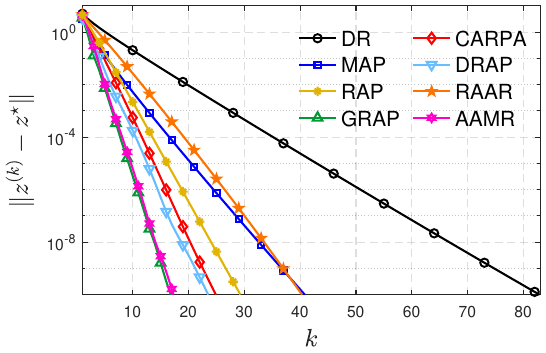}\label{fig:Experiment2-3}}
    \hskip3mm
    \subfloat[{$\tF=1.0$}]{\includegraphics[width=0.425\textwidth]{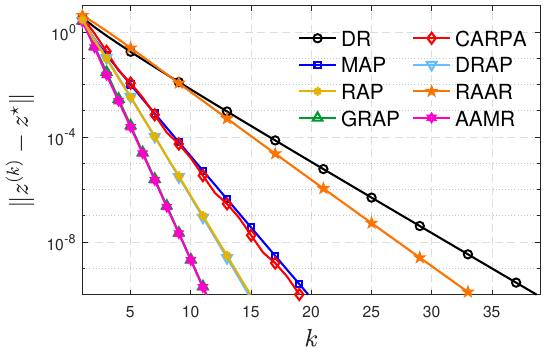}\label{fig:Experiment2-4}}
    \caption{Experimental results for different parameter sets.}
    \label{fig:Experiment2}
\end{figure}

The results are provided in Figure \ref{fig:Experiment2}, from which we observe that the practical performances of the methods verify Figure \ref{fig:W1comparenew}: 1) AAMR and GRAP are the fastest in all cases; 2) As the SP method is not compared, DR is the slowest one; 3) Our proposed method is comparable for $\tF=0.1,0.4,0.7$. 
For $\tF=1.0$, \carp is only marginally faster than MAP.
As the value of $\theta_F$ increases, the gap between these algorithms gradually narrows.
% \rc{not fully consistent with figure 2}

% \SynologyDrive\code\RPR-DR_twospace\Experiment2.m

\subsection{Feasibility of a ball tangent to a line}
%%%%%%%%%%%%%%%%%%%%%

In the first two-space polyhedral example, the performance of DR is not desirable. However, this is quite problem dependent, and in the next example we show that DR actually is the best. 
In this example, for the sake of simplicity, we constrain the problem in $\bbR^2$, and let $X$ be a line and $Y$ be a unit ball
\[
X = \Ba{x \mid  a^\top x = 1 } \qandq  Y = \Ba{x \mid \| x \|_{2} \leq 1 }
\]
where $a = (1, 1)^\top/\sqrt{2}$. Apparently, $X \cap Y = \Ba{ (1, 1)^\top/\sqrt{2} }$.

The following methods are considered for comparison: DR and nsDR, MAP and GRAP, the proposed \carp and its non-stationary version ns\carp. The settings of these algorithms are
\begin{itemize}
    \item The iteration of both DR and MAP are non-relaxed.
    
    \item For nsDR, we adopt the default setting as in~\eqref{eq:nsdr}. 

    \item The relaxation parameters GRAP are $\alpha_1=\alpha_2= 0.4$. 

    \item For \carp we set $\gamma = 0.5,\mu=1$, while for ns\carp (Algorithm \ref{alg:adcdr_adaptive1}) we use $\mu=1$, $(\gamma_0,\gamma_{\min},\gamma_{\max})=(0.5,0,1)$, $(c_1,c_2)=(0.5,50)$ and $\delta = 0.01$.
\end{itemize} 
The residual error $\norm{\zkp-\zk}$ is used, and four stopping criteria are considered $10^{\ba{-4,-6,-8,-10}}$. The maximum number of iterations for all methods is set as $10^4$, and the iteration is terminated when the tolerance is met or the maximum iteration is reached. 

As the solution of the problem $\xsol=(1, 1)^\top/\sqrt{2}$ is unique, we use $z^{(0)} = \xsol + 10 u ,~ \norm{u}=1$ as the starting point for all methods. To measure the average performance of the method, for each tolerance, we randomly generate $10^4$ starting points and compare the averaged number of iterations of each method. The result is provided in Table~\ref{tab:ball_line}, in which ``-'' means that the maximum iteration step is reached. From the table, we observe that: 1) Among the four stationary schemes, in contrast to Section \ref{otherprojection_methods}, DR is the fastest one; 2) The nsDR is the fastest over all methods; 3) ns\carp provides significant improvement over the standard \carp scheme, especially for tolerance $10^{-6,-8,-10}$, over an order acceleration is obtained. 
 
\begin{table}[!htb]
\centering  
% \begin{tabular}{|m{9mm}|m{9mm}|m{9mm}|m{16mm}|m{16mm}|m{16mm}|m{16mm}|}
\begin{tabular}{|c|c|c|c|c|c|c|}
\hline
 & DR  & nsDR & MAP & GRAP & \carp & ns\carp  \\ \hline
$10^{-4}$ & 24 & 15 & 292 & 178 & 104 & 64 \\ \hline
$10^{-6}$ & 177 & 21 & 6290 & 4481 & 3030 & 305 \\ \hline
$10^{-8}$ & 758 & 28 & 9995 & 9925 & 9172 & 790 \\ \hline
$10^{-10}$ & 1017 & 35 & - & 9989 & 9823 & 1140 \\ \hline
\end{tabular}
\caption{Comparison of number of steps required to reach stopping tolerance.}
\label{tab:ball_line}
\end{table}

%%%%%%%%%%%%%%%%%%%%%%%%%
\subsection{Sparse linear inverse problems}
The last example we consider is compressed sensing type problem, in which case the two sets are
\beq\label{eq:bp}
X = \Ba{ x\in \bbR^n \mid Ax = b  } \qandq Y = \Ba{x \in \bbR^n \mid  \|x\|_1 \leq c  }
\eeq
where $A \in \mathbb{R}^{m\times n}$ with full row-rank. Let $\xsol\in\bbR^n$ be a $\kappa$-sparse vector, we have and $b=A\xsol$ and $c=\norm{\xsol}_1$. 
The projection operator of $X$ is given by
$$
\proj_{X} (w) = w - A^\top(AA^\top)^{-1} (Aw-b),
$$
while the projection operator of $Y$ has a fast scheme as proposed in \cite{condat2016fast}.

\paragraph{A toy example}\label{sec:toy}
We first consider a toy example to compare the performance of the methods, for which we set $(m,n) = (500,2000)$ and $\kappa = 50$. The entries of $A$ are sampled from the normal distribution and normalized by $\sqrt{m}$. 
The methods DR/nsDR, \carp/ns\carp, MAP and GRAP are compared with the following settings
%methods in terms of fixed-point error $\|\zk-\zsol\|$ and the support of $\proj_{Y}(\zk)$:
\begin{itemize}
    \item The iterations of both DR and MAP are non-relaxed.
    
    \item For nsDR, we adopt the default setting as in~\eqref{eq:nsdr}. 

    \item The relaxation parameters GRAP are $\alpha_1=\alpha_2=0.75$. 

    \item For \carp we set $\gamma = 0.5,\mu=1$, while for ns\carp (Algorithm \ref{alg:adcdr_adaptive1}) we use $\mu=1$, $(\gamma_0,\gamma_{\min},\gamma_{\max})=(0.5,0,1)$, $(c_1,c_2)=(0.9,50)$ and $\delta = 0.01$.
\end{itemize}

 \begin{figure}[H]
    \centering
    \subfloat[{Convergence of $\norm{\zk-\zsol}$}]{\includegraphics[width=0.45\textwidth]{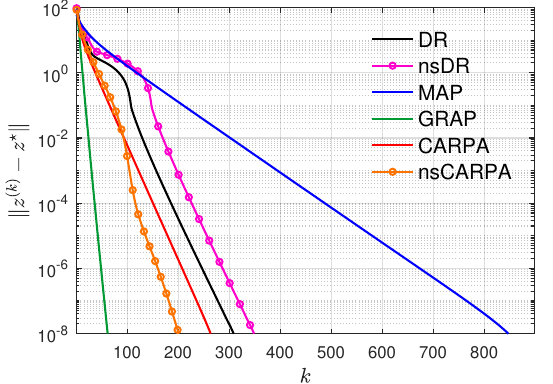} }
    \hskip3mm
    \subfloat[{The support of $\proj_{Y}(\zk)$}]{\includegraphics[width=0.45\textwidth]{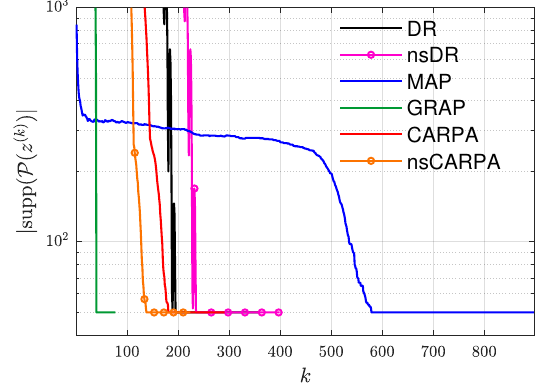} }
    \\[-3mm]
    \caption{Comparison on toy example: (a) Convergence of $\norm{\zk-\zsol}$; (b) Change of the support of $\proj_{Y}(\zk)$.} 
    \label{fig:toy-example}
\end{figure}
% \SynologyDrive\code\CMArP_basispursuit\cmp_l1normball_subspace.m

The results are presented in Figure \ref{fig:toy-example}. 
We observe that in Figure \ref{fig:toy-example}(a) GRAP exhibits the fastest performance as in the two-space feasibility problem. DR/nsDR and \carp/ns\carp are quite close to each other, with MAP being the slowest of all. 
In terms of the support size in Figure \ref{fig:toy-example}(b), the observations are close to the convergence of $\norm{\zk-\zsol}$, the main difference among the methods is the number of iterations needed for $\proj_{Y}(\zk)$ to find the support of $\xsol$.

\paragraph{Realistic examples}
We conclude this part by presenting experiments on some more realistic problems. 
Given a vector $f\in\bbR^n$ and a square matrix $B\in\bbR^{n\times n}$, suppose $f$ has a sparse representation under the columns of $B$, that is, there exists a sparse vector $\xsol$ such that $f = B \xsol$. Given a measurement matrix $M\in\bbR^{m\times n}$ and $b = M f = MB\xsol$, then recovering $f$ from $b$ can be obtained by solving \eqref{eq:bp}  with $A=MB$.
In this example, four problems cases listed in Table \ref{tab:cs-p} are considered, and we refer to \cite{van2009algorithm} for more details.

\begin{table}[htbp]
\centering
\begin{tabular}{|c|c|c|c|c|c|c| }
\hline
 problem & $(m,n)$ & $\kappa$ &   $B$ & $M$ \\ \hline
1& (1024,2048) & 120 & DCT & Dirac \\ \hline % 003
2& (600,2560) & 20 &  ${\rm Id}$ & Gaussian  \\ \hline % 007
3&(256,1024) & 32 & ${\rm Id}$ & Gaussian \\ \hline % 011
4&(200,1000) & 3 & DCT & oprestricted \\ \hline % 902
\end{tabular}
\caption{Problems cases and settings. \textit{DCT} means the discrete cosine transform. The term \textit{oprestriceted} represents the restriction operator, which extracts specific entries from a given vector. }% and returns a correspondingly shortened vector. }
\label{tab:cs-p}
\end{table}

\begin{figure}[H]
    \centering
    \subfloat[Problem 1]{\includegraphics[width=0.465\textwidth]{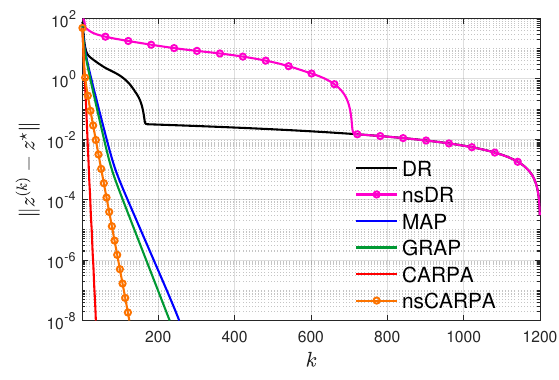}\label{fig:CS-1}}
    \hskip3mm
    \subfloat[Problem 2]{\includegraphics[width=0.465\textwidth]{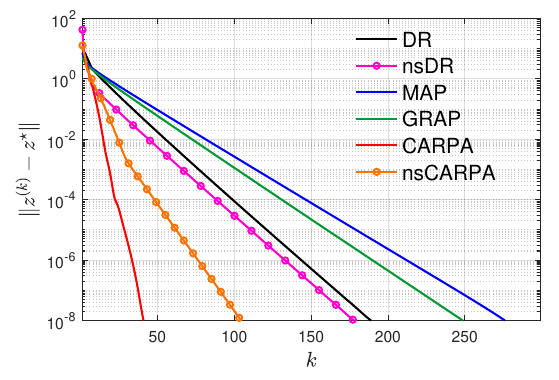}\label{fig:CS-2}}
    \\[-3mm]
    \subfloat[Problem 3]{\includegraphics[width=0.465\textwidth]{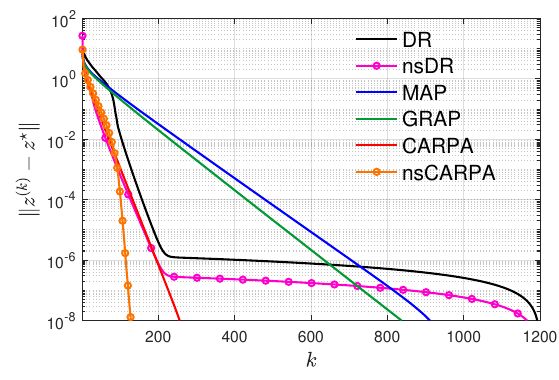}\label{fig:CS-3}}
    \subfloat[Problem 4]{\includegraphics[width=0.465\textwidth]{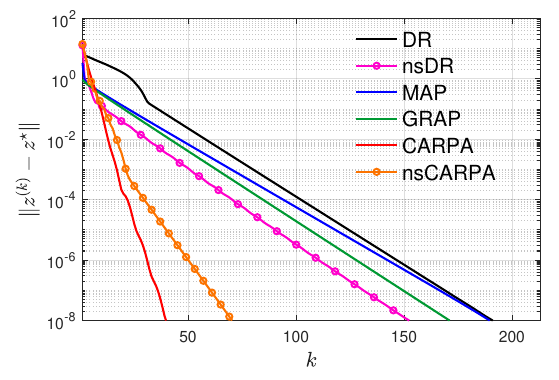}\label{fig:CS-4}}
    \caption{Experimental results for different problems.}
    \label{fig:Experiment4}
\end{figure}

We compare the performance of the six methods under the same parameter settings as outlined in the toy example. The observations in Figure \ref{fig:Experiment4} are quite different from those in Figure \ref{fig:toy-example}: 1) Overall \carp/ns\carp provide the best performance; 2) MAP and GRAP are similar to each other; 3) The performance of DR and nsDR is close, but unstable compared to other methods.

\section{Conclusions}
In this paper, we proposed a composed alternating relaxed projection algorithm for solving the feasibility problem, and a non-stationary version of the methods was developed. Theoretical analysis for both schemes is provided. 
For the feasibility of two spaces, we also analyzed the optimal convergence rate of the proposed method. 
Numerical experiments are provided to demonstrate the performance of the methods against several existing algorithms in the literature. Our proposed schemes provide comparative performances. 
In the future, we plan to extend the algorithm to the infeasible scenario and nonconvex setting.

\paragraph{Acknowledgements}
JL is supported by the National Natural Science Foundation of China (No. 12201405), the ``Fundamental Research Funds for the Central Universities'', the National Science Foundation of China (BC4190065) and the Shanghai Municipal Science and Technology Major Project (2021SHZDZX0102).

\begin{small}
\appendix
\section{Convergence rate of \carp for two subspaces}\label{appendix}

In this appendix, we provide the optimal convergence rate analysis for \carp. For RAAR and DRAP, we also provide a simple analysis at the end of the section. For the two subspaces $X$ and $Y$, we suppose that 
$\dim(X) = p$, $\dim(Y) = q$ such that $p+q< n$ and $p<q$ (for the case $p+q>n$ will be discussed later), then owing to \cite{bauschke2016optimal}, there exists an orthogonal matrix $D \in \mathbb{R}^{n \times n}$ such that 
\beq\label{eq:pX-pY}
\proj_X= D
\begin{bmatrix}
 \Id_{p} & 0 & 0 & 0 \\
0 & 0_p & 0 & 0 \\
0 & 0 & 0_{q-p} & 0 \\
0 & 0 & 0 & 0_{n-p-q}   
\end{bmatrix}
 D^* \qandq
\proj_Y = D
\begin{bmatrix}
C^2 & C S & 0 & 0 \\
C S & S^2 & 0 & 0 \\
0 & 0 & \Id_{q-p} & 0 \\
0 & 0 & 0 & 0_{n-p-q}
\end{bmatrix}
D^* ,
\eeq
where $C$ and $S$ are two $p \times p$ diagonal matrices defined by
$$
C:=\diag \Pa{ \cos (\theta_1), \ldots, \cos (\theta_p) } 
\qandq 
S:=\diag \Pa{ \sin (\theta_1), \ldots, \sin (\theta_p) }
$$
with the principal angles $\theta_1, \ldots, \theta_p$ between $X$ and $Y$ found in Definition \ref{def:principal angle}. For simplicity, we denote
\[
(c_i, s_i) = \Pa{ \cos (\theta_i), \sin(\theta_i) } ,~~ i=1,...,p .
\]
Moreover, we have 
\[
2\proj_X - \Id = \proj_X - \proj_{X^\bot}
\]
where $X^\bot$ is the orthogonal complement of $X$.

\subsection{Optimal convergence rate of \carp}
\label{pf:Fgammamu}

% \paragraph{Representation of the fixed-point operator}
From \eqref{eq:fcarp} we have
\[
\calF_{1} 
= \proj_Y\proj_X + \proj_{Y^\bot} \proj_{X^\bot} \\
\qandq 
\calF_{2} = (\proj_Y - \proj_{Y^\bot}) \proj_X , \\
\]
Consequently the fixed-point operator of \carp reads
\[
\begin{aligned}
 \fCDR 
 &= \proj_Y\proj_X + \proj_{Y^\bot} \proj_{X^\bot} - \gamma \proj_{X^\bot}  \\
 \qandq
 \fCDRmu &= (1-\mu) \Id  + \mu\Pa{ \proj_Y \proj_X+\proj_{Y^{\bot}} \proj_{X^{\bot}} - \gamma \proj_{X^{\bot}} } . 
\end{aligned}
\]
Combine with \eqref{eq:pX-pY}, we get
\begin{equation}
\begin{aligned}
\fCDRmu
% &= (1-\mu) \Id  + \mu(\proj_Y \proj_X+\proj_{Y^{\bot}} \proj_{X^{\bot}} - \gamma \proj_{X^{\bot}} )\\
& = D
\begin{bmatrix}
(1-\mu) \Id_p  + \mu C^2   &  - \mu C S   & 0 &0\\
\mu CS   &  (1-\mu -\mu \gamma ) \Id_p + \mu C^2 & 0 &0\\
0 & 0 & (1-\mu -\gamma\mu) \Id_{q-p} &0\\
0 & 0 & 0 &(1-\mu \gamma) \Id_{n-p-q} 
\end{bmatrix}
D^* . \\
% & = D {\rm blkdiag}(A, B, C) D^* ,
\end{aligned}
\end{equation}
The eigenvalues of $\fCDRmu$ is then determined by the block diagonal matrix. Note that when $\mu < 2/(1+\gamma)$, then $\gamma\mu< 1$, $-1< 1-\mu-\gamma\mu < 1$ and $1-\gamma\mu< 1$. As a result, we only need to study the eigenvalues of 
\begin{equation}\label{eq:mtx-M}
M = 
\begin{bmatrix}
    (1-\mu) \Id_p  + \mu C^2   &  - \mu C S  \\ 
    \mu CS   &  (1-\mu -\mu \gamma ) \Id_p  + \mu C^2
\end{bmatrix} .
\end{equation}
As $C, S$ are diagonal, studying the eigenvalues of $M$ in \eqref{eq:mtx-M} is equivalent to studying the eigenvalue of the following $2\times 2$ matrix, 
\[
N = \begin{bmatrix}
1-\mu + \mu c^2  & -\mu c s \\
  \mu c s & 1-\mu -\mu\gamma + \mu c^2 
\end{bmatrix} 
\]
where $c, s \in [0,1]$ satisfying $c^2+s^2=1$. 
It is trivial to see that $N$ has the following two eigenvalues
\[
\lambda_{1,2} = 1 - \mu + \mu c^2 +  \frac{ -\gamma\mu  \pm \mu \sqrt{ \gamma^2 - 4 c^2s^2 }  }{2}  . 
\]
We need to study the dependencies of $\lambda$ over $\gamma, \mu$, for which we have the following proposition. 

From now on, we denote intervals $$\calI_1 = \Big[0, \sfrac{1-\sqrt{1-\gamma^2}}{2}\Big], ~~~ \calI_2 = \Big[\sfrac{1+\sqrt{1-\gamma^2}}{2},1\Big]\qandq \calI_3 = [0,1]\verb|\| (\calI_1\cup \calI_2).$$
Other denotation includes
$$
\stl = \sfrac{1+\sqrt{1-\gamma^2}}{2}, ~~~\mu_0 = \sfrac{2}{1+\gamma +\sqrt{1-\gamma^2}}, ~~~\varphi= \max\{|1-\mu\gamma|,|1-\mu-\mu\gamma|\}.
$$

\begin{proposition}
\label{th:Fgammamu}
 Let $\gamma \in [0,1]$, $\mu\in ]0,2/(1+\gamma)[$, and  $p+q\leq n$. Then, the \carp mapping  $\fCDRmu=(1-\mu) \Id  + \mu(\proj_Y \proj_X+\proj_{Y^{\bot}} \proj_{X^{\bot}} - \gamma \proj_{X^{\bot}} )$ is linearly convergent to $ \proj_{ {\Fix} \fCDR } = \proj_{X\cap Y}$. Besides, when $\mu=1$, $p=q$ and $p+q = n$, the optimal convergence rate would be
 $$
 \rsol= \left\{
\begin{array}{ll}
\cF^2-\cF\sF, & \qifq \sF^2 \leq \tfrac{1+\ssqrt{1-\gamma^2}}{2},\\
 c_ps_p-c_p^2 , & \qifq \sF^2 > \tfrac{1+\ssqrt{1-\gamma^2}}{2}
\end{array}
\right.
 $$
 with 
 $$
 \gmsol= \left\{
\begin{array}{ll}
2\cF\sF, & \qifq \sF^2 \leq \tfrac{1+\ssqrt{1-\gamma^2}}{2},\\
 2c_ps_p, & \qifq \sF^2 > \tfrac{1+\ssqrt{1-\gamma^2}}{2}.
\end{array}
\right.
 $$
\end{proposition}
 
\begin{proof}
Depends on the sign of $\gamma^2 - 4 c^2s^2$, we have the following two cases
\begin{itemize}
\item {\bf Case $\gamma^2 - 4 c^2s^2 < 0$:} 
For this case, we have
\[
\gamma^2 - 4 c^2s^2
= \gamma^2 - 4 s^2 (1-s^2) 
= 4s^4 - 4s^2 + \gamma^2
< 0 
\]
which means 
\[
s^2 \in \calI_3 \subset [0, 1] . 
\]
% \Big[\sfrac{1-\sqrt{1-\gamma^2}}{2},~ \sfrac{1+\sqrt{1-\gamma^2}}{2}\Big]
For the eigenvalues, we have the following simplified expression of the magnitude of the eigenvalues
\[
\begin{aligned}
\abs{\lambda_{1,2}}^2
&=   (1+\gamma) s^2 \mu^2 -  (2s^2+\gamma) \mu + 1 \\
&=  \Pa{  (1+\gamma)  \mu - 2 } \mu s^2 -  \gamma \mu + 1 \\
\end{aligned}
\]
Note that as $\mu<2/(1+\gamma)$, we have $(\gamma+1) \mu - 2 < 0$ and 
$$
\begin{aligned}
\ell(s^2) 
&= \Pa{  (1+\gamma)  \mu - 2 } \mu s^2 -  \gamma \mu + 1
\end{aligned}
$$
is decreasing for $s^2 \in [0,1]$. Hence we only need to consider $s^2=\frac{1-\sqrt{1-\gamma^2}}{2}$ and $s^2=\frac{1+\sqrt{1-\gamma^2}}{2}$, and we have
\begin{itemize}
    \item $s^2=\frac{1-\sqrt{1-\gamma^2}}{2}$: again $2-\mu(1+\gamma) > 0$ and
        \[
        \begin{aligned}
        \ell\Pa{ \tfrac{1-\sqrt{1-\gamma^2}}{2}} 
        &= \Pa{  (1+\gamma)  \mu - 2 } \mu \sfrac{1-\sqrt{1-\gamma^2}}{2} -  \gamma \mu + 1   \\
        &= \sfrac{(1+\gamma)(1-\sqrt{1-\gamma^2})}{2} \mu^2 - \pa{1+\gamma-\sqrt{1-\gamma^2}} \mu + 1 ,
        \end{aligned}
        \]
        which is quadratic function of $\mu$ and denoted as $q(\mu)$. 
        Note that 
        \[
        \begin{aligned}
            \sfrac{ 1+\gamma-\sqrt{1-\gamma^2} }{(1+\gamma)(1-\sqrt{1-\gamma^2})} - \sfrac{2}{1+\gamma}
            &= \sfrac{ 1+\gamma-\sqrt{1-\gamma^2} - 2+ 2\sqrt{1-\gamma^2} }{(1+\gamma)(1-\sqrt{1-\gamma^2})}  \\
            &= \sfrac{ \gamma - 1 + \sqrt{1-\gamma^2} }{(1+\gamma)(1-\sqrt{1-\gamma^2})}  > 0 .
        \end{aligned}
        \]
        This means $q(\mu)$ is decreasing for $\mu \in ]0, 2/(1+\gamma)[$. 
        Since $q(0) = 1$ and 
        \[
        \begin{aligned}
            q\Pa{2/(1+\gamma)}
            &= \Pa{  (1+\gamma)  \mu - 2 } \mu \sfrac{1-\sqrt{1-\gamma^2}}{2} -  \gamma \mu + 1   \\
           &= 1 - \sfrac{2 \gamma}{1+\gamma}
           = \sfrac{1 - \gamma}{1+\gamma}
            < 1 .
        \end{aligned}
        \]
        We get $\ell(s^2) < 1$ as long as $\mu \in ]0, 2/(1+\gamma)]$. 

    \item $s^2=\frac{1+\sqrt{1-\gamma^2}}{2}$: follow the discussions above, we have
        \[
        \begin{aligned}
        \ell\Pa{ \tfrac{1+\sqrt{1-\gamma^2}}{2}} 
        &= \Pa{  (1+\gamma)  \mu - 2 } \mu \sfrac{1+\sqrt{1-\gamma^2}}{2} -  \gamma \mu + 1   \\
        &= \sfrac{(1+\gamma)(1+\sqrt{1-\gamma^2})}{2} \mu^2 - \pa{1+\gamma+\sqrt{1-\gamma^2}} \mu + 1 ,
        \end{aligned}
        \]
        which again is quadratic function of $\mu$ and denoted by $q(\mu)$. 
        Note that 
        \[
        \begin{aligned}
            \sfrac{ 1+\gamma+\sqrt{1-\gamma^2} }{(1+\gamma)(1+\sqrt{1-\gamma^2})} - \sfrac{2}{1+\gamma}
            &= \sfrac{ 1+\gamma+\sqrt{1-\gamma^2} - 2 - 2\sqrt{1-\gamma^2} }{(1+\gamma)(1+\sqrt{1-\gamma^2})}  \\
            &= \sfrac{ \gamma - 1 - \sqrt{1-\gamma^2} }{(1+\gamma)(1+\sqrt{1-\gamma^2})}  < 0 .
        \end{aligned}
        \]
        This means $q(\mu)$ is decreasing for $\mu \in \big]0, \frac{ 1+\gamma+\sqrt{1-\gamma^2} }{(1+\gamma)(1+\sqrt{1-\gamma^2})} \big[$ and increasing for $\mu \in \big[ \frac{ 1+\gamma+\sqrt{1-\gamma^2} }{(1+\gamma)(1+\sqrt{1-\gamma^2})}, 2/(1+\gamma) \big]$. 
        Since $q(0) = 1$ and $q\Pa{2/(1+\gamma)} = \frac{1 - \gamma}{1+\gamma}$, then 
        \[
        \begin{aligned}
            q\bPa{ \sfrac{ 1+\gamma+\sqrt{1-\gamma^2} }{(1+\gamma)(1+\sqrt{1-\gamma^2})} }
            &= \sfrac{(1+\gamma)(1+\sqrt{1-\gamma^2})}{2} \mu^2 - \pa{1+\gamma+\sqrt{1-\gamma^2}} \mu + 1  \\
            &= 1 - \sfrac{ (1+\gamma+\sqrt{1-\gamma^2})^2 }{2(1+\gamma)(1+\sqrt{1-\gamma^2})}  = 0  .
        \end{aligned}
        \]
        Therefore, we conclude that $\ell(s^2) < 1$ as long as $\mu \in ]0, 2/(1+\gamma)]$. 
        
\end{itemize}

\item {\bf Case $\gamma^2 - 4 c^2s^2 \geq 0$:}
For this case, we have
\[
\gamma^2 - 4 c^2s^2
= \gamma^2 - 4 s^2 (1-s^2) 
= 4s^4 - 4s^2 + \gamma^2
> 0 
\]
which means 
\[
s^2 \in \calI_1\cup \calI_2 \subset [0, 1] . 
\]
For the eigenvalues, we have

\[
\begin{aligned}
\lambda_{1,2}
&= 1 - \mu + \mu c^2 +  \frac{ -\gamma\mu  \pm \mu \sqrt{ \gamma^2 - 4 c^2s^2 }  }{2}  \\
% &= 1 + \mu\Pa{ - 1 +  c^2 +  \frac{ -\gamma  \pm  \sqrt{ \gamma^2 - 4 c^2s^2 }  }{2}  } \\
&= 1 - \mu \frac{ \gamma+2s^2  \pm  \sqrt{ \gamma^2 - 4s^2 + 4s^4 }  }{2}  .
\end{aligned}
\]
We have the following cases, %let $\bar{s} = s^2$
\begin{itemize}
    \item For $f_1(s^2) = \gamma+2s^2  -  \sqrt{ \gamma^2 - 4s^2 + 4s^4 } $. Note that 
    \[
    (\gamma+2s^2)^2 - (\gamma^2 - 4s^2 + 4s^4) = 4(1+\gamma) s^2 > 0 . 
    \]
    Hence $f_1(s^2)> 0$. Note that $f_1(s^2)$ is increasing for $s^2 \in \calI_1$ and decreasing in $s^2 \in \calI_2$, thus achieves maximum value for $s^2=\frac{1+\sqrt{1-\gamma^2}}{2}$ with 
    \[
      f_1(\frac{1+\sqrt{1-\gamma^2}}{2}) = 1+\sqrt{1-\gamma^2} +\gamma > 2.
    \]
    However
    \[
    1-\mu\sfrac{f_1(\frac{1+\sqrt{1-\gamma^2}}{2})}{2} = 1-\frac{\mu(1+\gamma +\sqrt{1-\gamma^2})}{2} >-1
    \]
    because $\mu(1+\gamma)<2.$
    % Note that $f_1(s^2)$ may be negative however 
    % \[
    % f_1(1) = 2 .
    % \]
    % Then we have 
    In addition,
    \[
    1-\mu\sfrac{f_1(2)}{2} = 1 - \mu \in ]-1, 1[.
    \]
    Hence $\lambda_1 = 1-\frac{\mu f_1(s^2)}{2} \in ]-1,1[$.

      \item Now we define $f_2(s^2) = \gamma+2s^2  +  \sqrt{ \gamma^2 - 4s^2 + 4s^4 } $, it can be verified that $f_2(s^2)$ is decreasing for $s^2\in \calI_1$ and increasing for $s^2 \in [1/2, 1]$, hence achieves maximum value for $s^2=1$, for which we have
    \[
    f_2(1) = 2(1+\gamma) > f_2(0) = 2\gamma.
    \]
    Then for the smaller eigenvalue we have
    \[
    1-\mu\sfrac{f_2(1)}{2}  
    = 1 - \mu \frac{ 2(1+\gamma) }{2}
    = 1 - \mu (1+\gamma)
    > 1-2 = -1 .
    \]
    Hence $\lambda_2 = 1-\frac{\mu f_2(s^2)}{2} \in ]-1,1[$.
\end{itemize}

\end{itemize}

To study square of the modulus of eigenvalues, we consider the function $g_{\gamma, \mu,r}: [0,1] \rightarrow \mathbb{R}$ given by 
\beq \label{gfcn}
g_{\gamma,\mu,r}(s^2): =\begin{cases} 
\Pa{  (1+\gamma)  \mu - 2 } \mu s^2 -  \gamma \mu + 1,&\qifq  s^2 \in \calI_3,\\
 (1 - \mu \frac{ \gamma+2s^2  +(-1)^r  \sqrt{ \gamma^2 - 4s^2 + 4s^4 }  }{2})^2,& \qifq s^2\in \calI_1\cup \calI_2 .\end{cases}
\eeq

Discussion about the monotonicity of  $g_{\gamma,\mu,r}(s^2), r=1,2$:
\begin{itemize}
    \item $g_{\gamma,\mu,1}(s^2)$ is always decreasing on $\calI_1\cup \calI_3$, however for interval $\calI_2$ it is increasing when $0<\mu<\mu_0$,  decreasing when $1< \mu <\frac{2}{1+\gamma}$ and first decreasing then increasing when $\mu_0 < \mu <1$ with a zero in $\calI_2$;
    % \rc{(it is because $1-\mu\sfrac{f(s^2)}{2}$)}
    \item $g_{\gamma,\mu, 2}(s^2)$ is always increasing on $\calI_1$ and decreasing on $\calI_3$, however for interval $\calI_2$ it is decreasing when $0<\mu<\frac{1}{1+\gamma}$,  increasing when $\mu_0=\frac{2}{1+\gamma+\sqrt{1-\gamma^2}}< \mu <\frac{2}{1+\gamma}$ and first decreasing then increasing when $\frac{1}{1+\gamma} < \mu <\mu_0$ with a zero in $\calI_2$.
\end{itemize}

Denote $g_{\gamma,\mu} =\max\{g_{\gamma,\mu,1},g_{\gamma,\mu,2}\}$ and what we want is  $g_{\max} =\max_s g_{\gamma,\mu}(s^2)$ which means $|\lmax(M)| = \sqrt{g_{\max}}$
\begin{itemize}
    \item for all  $0\leq \us \leq \os  \leq \stl  $, 
\begin{equation}\label{mono1}
    (1-\mu\gamma)^2 = g_{\gamma,\mu,1}(0) \geq g_{\gamma,\mu,1}(\us)\geq g_{\gamma,\mu,1}(\os)\geq  g_{\gamma,\mu,2}(\stl).
\end{equation} 
\item for all $ \stl \leq  \us \leq\os  \leq 1$,
\begin{equation}\label{mono2}
 (1-\mu)^2 =g_{\gamma,\mu,1}(1) \geq g_{\gamma,\mu,1}(\os)\geq g_{\gamma,\mu,1}(\us)\geq  g_{\gamma,\mu,1}(\stl) , \qifq \mu \le\mu_0,
\end{equation}
and
\begin{equation}\label{mono3}
 (1-\mu(1+\gamma))^2=g_{\gamma,\mu,2}(1)\geq g_{\gamma,\mu,2}(\os)\geq g_{\gamma,\mu,2}(\us)\geq  g_{\gamma,\mu,2}(\stl), \qifq\mu >\mu_0. 
\end{equation}
\end{itemize}
 
\begin{figure}[ht]
    \centering
    \includegraphics[width=0.4\linewidth]{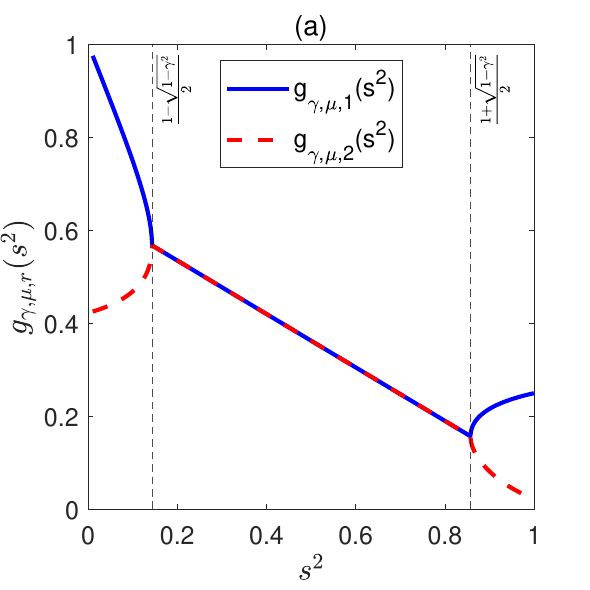}
    \includegraphics[width=0.4\linewidth]{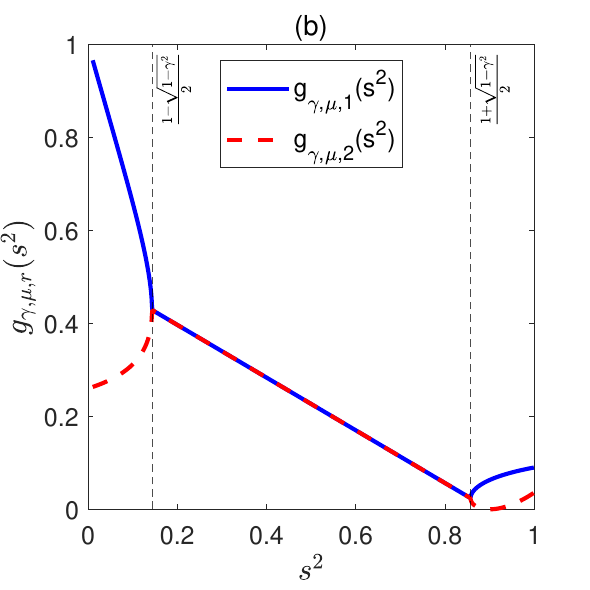}\\
     \includegraphics[width=0.4\linewidth]{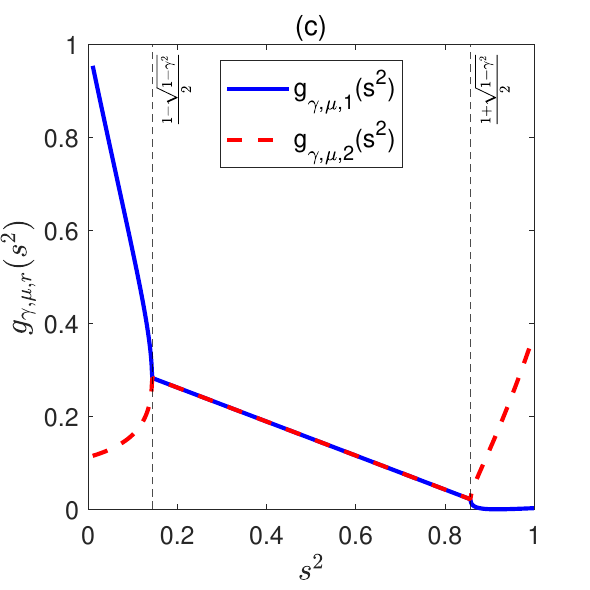}
    \includegraphics[width=0.4\linewidth]{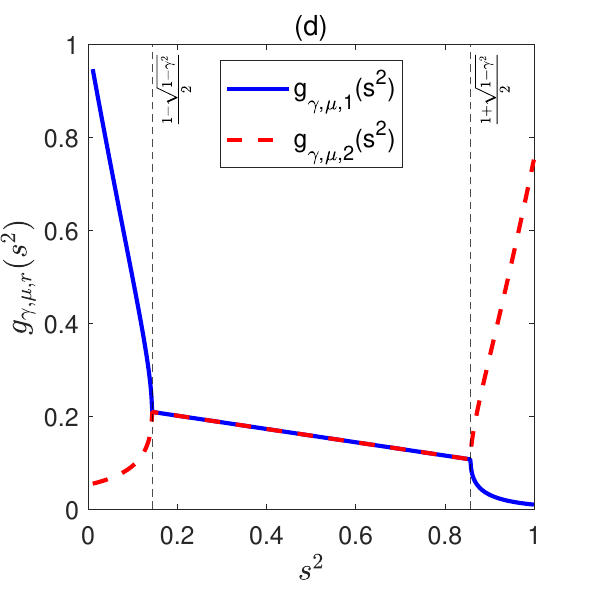}
    \caption{Illustration of $g_{\gamma,\mu,r}(s^2)$:
    (a) $0<\mu< \frac{1}{1+\gamma}$, (b)$\frac{1}{1+\gamma}<\mu <\mu_0$,(c)$\mu_0<\mu< 1$, (d) $1<\mu< \frac{2}{1+\gamma}$.}
    \label{mono}
\end{figure}

Deduce from the monotonicity \eqref{mono1}, \eqref{mono2} and \eqref{mono3}, we have the eigenvalue with the largest absolute value of M :
\begin{itemize}
    \item $\sF^2\in \calI_1, s_p^2\in \calI_2$ 
    \begin{itemize}
        \item when $\mu<\mu_0$  
            \beqn
                |\lmax|
                = 1 - \mu \frac{ \gamma+2\sF^2  -  \sqrt{ \gamma^2 - 4\sF^2 + 4\sF^4 }  }{2}.
            \eeqn
        \item  when $\mu>\mu_0 $ 
            \beqn
                |\lmax|= \max  \{ 1 - \mu \frac{ \gamma+2\sF^2  -  \sqrt{ \gamma^2 - 4\sF^2 + 4\sF^4 }  }{2}, -1 + \mu \frac{ \gamma+2s_p^2 +\sqrt{ \gamma^2 - 4s_p^2 + 4s_p^4 }  }{2} \}.
            \eeqn
    \end{itemize}
     \item $\sF^2\in \calI_1, s_p^2\in \calI_1\cup \calI_3$ 
           \beqn
                |\lmax|
                =1 -\mu\frac{\gamma+2\sF^2-\sqrt{ \gamma^2 - 4\sF^2 + 4\sF^4 }  }{2}.
            \eeqn
    \item $\sF^2\in \calI_3, s_p^2\in \calI_3$
            \beqn
                |\lmax|
                =\sqrt{\Pa{(1+\gamma)\mu - 2 } \mu\sF^2-\gamma\mu +1}.
            \eeqn
    \item $\sF^2\in \calI_3, s_p^2\in \calI_2$
         \begin{itemize}
         \item when $\mu<\mu_0$  
              \beqn
                |\lmax|= \max\{\sqrt{\Pa{(1+\gamma)\mu-2} \mu \sF^2-\gamma\mu +1}, 1-\mu\frac{ \gamma+2 s_p^2 -\sqrt{ \gamma^2 - 4 s_p^2 + 4 s_p^4 }  }{2}\}.
            \eeqn
         \item  when $\mu>\mu_0 $ 
              \beqn
                |\lmax|= \max\{\sqrt{\Pa{(1+\gamma)\mu-2} \mu \sF^2-\gamma\mu +1}, -1+\mu\frac{ \gamma+2 s_p^2 +\sqrt{ \gamma^2 - 4 s_p^2 + 4 s_p^4 }  }{2} \}.
            \eeqn
        \end{itemize}
    \item $\sF^2\in \calI_2, s_p^2\in \calI_2$
     \begin{itemize}
        \item when $\mu<\mu_0$  
            \beqn
                |\lmax|
                = \max  \{ 1 - \mu \frac{ \gamma+2\sF^2  -  \sqrt{ \gamma^2 - 4\sF^2 + 4\sF^4 }  }{2}, 1 - \mu \frac{ \gamma+2s_p^2 -\sqrt{ \gamma^2 - 4s_p^2 + 4s_p^4 }  }{2} \}.
            \eeqn
        \item  when $\mu>\mu_0 $ 
            \beqn
                |\lmax|= \max  \{ 1 - \mu \frac{ \gamma+2\sF^2+\sqrt{ \gamma^2 - 4\sF^2 + 4\sF^4 }  }{2}, -1+\mu \frac{ \gamma+2s_p^2 +\sqrt{ \gamma^2 - 4s_p^2 + 4s_p^4 }  }{2} \}.
            \eeqn
    \end{itemize}
\end{itemize}

Now we think about fixed $\mu\equiv 1$ to discover the optimal $\gmsol$ and corresponding optimal rate $\rsol$ when $p=q$ and $p+q=n$. In this case, it must be (c) in Figure \ref{mono} as $\frac{1}{1+\gamma}<\mu_0<\mu$.

When $0<\sF^2< s_p^2 < \sfrac{1+\sqrt{1-\gamma^2}}{2}$, $\xi(\fCDR)=\lambda_{1}(\sF^2)$.
Minimizing $\lambda_{1}(\sF^2)$ yields $\gmsol = 2\cF \sF$, with optimal convergence rate is $\rsol = \cF ^2-\cF \sF $. 

When $\sF^2 < \sfrac{1+\sqrt{1-\gamma^2}}{2} < s_p^2$, $\xi(\fCDR)=\max\{\lambda_1(s^2_F), -\lambda_2(s^2_p)\}$. The optimal convergence rate is still $\rsol = \cF ^2-\cF \sF $ when $\gmsol = 2\cF \sF$ because minimizing $\lambda_{1}(\sF^2)$ leads to the boundary point of $\calI_1$ and $\calI_2$, minimizing $-\lambda_{2}(\sF^2)$ leads to the boundary point of $\calI_2$ and $\calI_3$ and it is decreasing on $\calI_2$.

When $\sfrac{1+\sqrt{1-\gamma^2}}{2}  < \sF^2 < s_p^2$ , $\xi(\fCDR)=-\lambda_2(s_p^2)$.
Minimizing $-\lambda_2(s_p^2)$ over $\gamma$ yields$\gmsol=2c_ps_p$, with optimal convergence rate is $\rsol = c_ps_p -c_p^2$.  
\end{proof}
\begin{remark}
If $ p+q \geq n$. We may find some $i \in \bbN$ such that $n^{\prime}:= n+i>p+q$. Define $X^\prime:= X \times\left\{0_i\right\} \subset \bbR^{n^\prime}, Y^\prime:=Y \times\left\{0_i\right\} \subset \bbR^{n^\prime}$, and $\fCDRmu^\prime   = (1-\mu)\Id + \mu( \proj_{Y^\prime} \proj_{X^\prime}  +\proj_{Y^{\prime\bot}} \proj_{X^{\prime\bot}} - \gamma \proj_{X^{\prime\bot}}). $
It is clear that $1 \leq p=\dim (X^{\prime}) \leq \dim (Y^{\prime})=q$ and $p+q<n^{\prime}$. Observe from Definition \ref{def:principal angle} that the principal angles between $X^{\prime}$ and $Y^{\prime}$ are the same with the ones between $X$ and $Y$. Moreover, we have $\proj_{X^{\prime}}=
\begin{bmatrix}\proj_X & 0 \\ 0 & 0_k
\end{bmatrix}, 
\proj_{Y^{\prime}}=
\begin{bmatrix}\proj_Y & 0 \\ 0 & 0_k\end{bmatrix}$,
$\proj_{X^{\prime \bot} }  = \Id - \proj_{X^{\prime}} = \begin{bmatrix}\proj_{X^\bot} & 0 \\ 0 & \Id\end{bmatrix},
 \proj_{Y^{\prime \bot} }  = \Id - \proj_{Y^{\prime}} = 
 \begin{bmatrix}\proj_{Y^\bot} & 0 \\ 0 & \Id
 \end{bmatrix}$
and thus
\begin{equation}
  \fCDRmu^{\prime}=\begin{bmatrix}
\fCDRmu & 0 \\
0 & (1-\gamma) \Id
\end{bmatrix}. 
\label{fcdrmu}
\end{equation}
From Proposition \ref{th:Fgammamu}  that $\fCDRmu^\prime $ is normal, and so is $\fCDRmu$. From formula \eqref{fcdrmu}, we can get that $\fCDRmu$ is convergent if and only if $\fCDRmu^\prime $ is convergent with the same rate. 
\end{remark}

\subsection{Optimal convergence rates for DRAP and RAAR}
\label{drapraar}

In this part, we provide a simple derivation of the optimal convergence rates of DRAP and RAAR in two-subspace feasibility problem burt only consider when $p = q = \frac{n}{2}$. 

\paragraph{{A relaxed Douglas--Rachford (DRAP) \cite{thao2018convergent}}} DRAP is another convex relaxation of DR, which can be viewed as a convex combination of DR and MAP, with $\mu \in (0,1]$. It is defined as 
 $$
 \calF_{_{\rm DRAP}} = \mu(\proj_Y \proj_X+\proj_{Y^{\bot}} \proj_{X^{\bot}}) + (1-\mu)\proj_Y \proj_X = \proj_Y \proj_X + \mu \proj_{Y^{\bot}} \proj_{X^{\bot}}
 $$
 and can also be expressed in matrix form as
 \begin{equation}
 \calF_{_{\rm DRAP}}= D
 \begin{bmatrix}
 C^2   & -\mu C S \\
 CS  & \mu C^2 \\ 
\end{bmatrix}
D^*.
\end{equation}

We get from the above format and the block determinant formula that 
$$
\sigma\left(\calF_{_{\rm DRAP}}\right)=\left\{ \lambda_{k_1},\lambda_{k_2}     \mid k=1, \ldots, p\right\}
$$
where
$$\lambda_{k_r} =  \frac{(1+\mu) c^2_k}{2} +(-1)^r \frac{\sqrt{(1+\mu)^2c_k^4-4\mu c_k^2 }}{2}, \ r = 1,2, \quad k = 1,\cdots,p$$
and 
$$
 |\lambda_{k_r} | = \left\{
\begin{array}{ll}
\sqrt{\mu}c_k  , & \qifq \frac{(1-s_k)^2}{c_k^2} \leq \mu \leq1,\\
    | \frac{1+\mu}{2}c_k^2 + (-1)^r \frac{c_k}{2}\sqrt{(\mu+1)^2c_k^4-4\mu c_k^2}| , & \qelse.
\end{array}
\right.
$$
Denote 
\beq \label{xi}
\xi_k = \max\{\left| \lambda_{k_1} \right|,\left| \lambda_{k_2} \right|\},
\eeq
and to find optimal $\msol$ corresponding to the optimal convergence rate is equivalent to solve this optimization problem:
\beq \label{drap}
\min_{\mu} \max_k \{\xi_k(\mu)\mid k= 1,\cdots, p \}=  \min_{\mu} \xi_F(\mu).  
\eeq
The minimizer of \eqref{drap} is $\msol = \frac{(1-\sF)^2}{\cF^2} $ and the reason for the first equivalence is the monotonicity of $\xi_k$ w.r.t $k$. The corresponding optimal rate is 
$$
\xi_F(\msol) = \sqrt{\frac{(1-\sF)^2}{\cF^2}}\cF =1-\sF
$$
% where $\argmin_\mu \ \xi(\mu) = \frac{(1-s)^2}{c^2}$ and $\min_{\mu} \xi(\mu) = \sqrt{\frac{(1-s)^2}{c^2}}c =1-s$.
 \end{small}

\paragraph{Relaxed Averaged Alternating Reflections (RAAR) \cite{luke2004relaxed}} the relaxed averaged alternating reflections is proven to be convergent when $\mu \in (0,1]$. It is defined as 
\begin{equation}
    \calF_{_{\rm RAAR}} = \mu(\proj_Y \proj_X+\proj_{Y^{\bot}} \proj_{X^{\bot}}) + (1-\mu)\proj_X = D
\begin{bmatrix}
 \Id -\mu S^2   & -\mu C S\\
\mu CS  & \mu C^2 \\ 
\end{bmatrix} 
D^*.
\end{equation}

We get from the above format and the block determinant formula that 
$$
\sigma ( \calF_{_{\rm RAAR}} ) = 
\left \{  \lambda_{k_1}, \lambda_{k_2}\mid k=1, \ldots, p\right\}
$$
where
$$\lambda_{k_r} =  \frac{1-\mu+2\mu c^2_k}{2} +(-1)^r \frac{\sqrt{(\mu-1)^2-4\mu^2 c_k^2s_k^2}}{2}, \ r = 1,2, \quad k = 1,\cdots,p$$
and
$$
\left| \lambda_{k_r} \right| = \left\{
\begin{array}{ll}
    \sqrt{\mu}c_k, & \qifq \frac{1}{1+2c_ks_k}  < \mu \leq 1, \\
   | \frac{1-\mu +2\mu c_k^2}{2} +(-1)^r \frac{1}{2}\sqrt{(\mu-1)^2-4\mu^2 c_k^2s_k^2}|, & \qifq  0  < \mu \leq \frac{1}{1+2c_ks_k}.
\end{array}
\right.
$$
 It is because $\xi_k$ with same definition in \eqref{xi} is not monotone w.r.t $k$ when $0\leq\tF \leq \theta_p \leq \pi/2 $. We omit the complicated analysis here which is similar to the analysis of \carp. At least, we derive the optimal $\msol = \frac{1}{1+2c_ks_k} $  when $0\leq\tF \leq \theta_p \leq \pi/4 $ and the corresponding convergence rate is $\frac{\cF}{\sqrt{1+2\cF\sF}}$.

\small
\bibliographystyle{plain}
\bibliography{bib1}

\end{document}